\RequirePackage[leqno]{amsmath} 
\documentclass[leqno]{amsart}

\usepackage[utf8]{inputenc}     % Accept different input encodings
\usepackage[T1]{fontenc}

\usepackage{amsmath,amssymb,amsfonts,mathtools}
\usepackage{latexsym}
\usepackage{pifont}

\usepackage[margin=0.5in]{geometry}

\usepackage[inline]{enumitem} 
\setlist[itemize]{topsep=0pt,after=\vspace{1.5\baselineskip}}

\usepackage{color}
\usepackage[demo]{graphicx}
\usepackage{efbox}
\efboxsetup{linecolor=gray,linewidth=0.8pt}
\usepackage{floatrow}
\usepackage{float}
\usepackage{caption}
\usepackage{subcaption}
\usepackage[leftcaption]{sidecap}
\sidecaptionvpos{figure}{m}

\usepackage{empheq}
\usepackage{tabularx,multirow,array} 

\usepackage{multicol}
\usepackage{xparse}

\usepackage[bookmarksopen,colorlinks=true,pdfencoding=auto, psdextra]{hyperref}
\makeatletter

\ExplSyntaxOn
\int_new:N \g_tohi_const_int
\int_new:N \g_tohi_const_sub_int
\tl_new:N  \g_tohi_const_char_tl

\cs_new_protected:Nn \tohi_print_constant:nn 
{
	#1 \textsubscript {#2} 
} 

\NewDocumentCommand\resetconstants{m}
{
	\int_gincr:N \g_tohi_const_int
	\int_gzero:N \g_tohi_const_sub_int
	\tl_gset:Nn  \g_tohi_const_char_tl {#1}
}

\NewDocumentCommand\const{m}
{
	\tl_if_exist:cTF
	{
		c_tohi_const_\int_use:N\g_tohi_const_int _#1_tl
	}
	{
		\tl_use:c {c_tohi_const_\int_use:N\g_tohi_const_int _#1_tl }
	}
	{
		\int_gincr:N \g_tohi_const_sub_int
		\tl_const:cx {c_tohi_const_\int_use:N\g_tohi_const_int _#1_tl }
		{ \exp_not:N\tohi_print_constant:nn {\g_tohi_const_char_tl }{\int_use:N \g_tohi_const_sub_int}}
		\tl_use:c {c_tohi_const_\int_use:N\g_tohi_const_int _#1_tl }
	}
}

\ExplSyntaxOff
\newcommand{\inlineitem}[1][]{%
	\ifnum\enit@type=\tw@
	{\descriptionlabel{#1}}
	\hspace{\labelsep}%
	\else
	\ifnum\enit@type=\z@
	\refstepcounter{\@listctr}\fi       
	\quad\@itemlabel\hspace{\labelsep}%
	\fi}
\DeclarePairedDelimiter\abs{\lvert}{\rvert}
\DeclarePairedDelimiter\norm{\lVert}{\rVert}
\DeclarePairedDelimiter\tonda{(}{)}
\DeclarePairedDelimiter\quadra{[}{]}
\DeclarePairedDelimiter\graffa{\{}{\}}
\newcommand{\into}{\int_\Omega}

\setlist[itemize]{noitemsep, topsep=0pt}

\def\R{\mathbb R} \def\N{\mathbb N}

\def\R{\mathbb R} \def\N{\mathbb N} 
\def\TM{T_{max}} 
\def
\@cite
#1#2{[#1\if@tempswa, #2\fi]}
\makeatother

\newtheorem{theorem}{Theorem}[section]

\newtheorem{assumptions}[theorem]{Assumptions}

\newtheorem{lemma}[theorem]{Lemma}

\newtheorem{remark}{Remark}

\title[Uniform-in-time boundedness in a class of chemotaxis models]{Uniform-in-time boundedness in a class of local and nonlocal nonlinear attraction-repulsion chemotaxis models with logistics}

\author[Alessandro Columbu, Rafael Diaz Fuentes, Silvia Frassu]{Alessandro Columbu, Rafael Diaz Fuentes, Silvia Frassu$^\star$}

\subjclass[2020]{Primary: 35K55, 35Q92. Secondary:  92C17.}
\keywords{Chemotaxis, Attraction-repulsion, Nonlinear production, Boundedness. \\
	\textit{$^\star$Corresponding author}: silvia.frassu@unica.it}

\begin{document}
	\maketitle
	{
		\medskip
		\centerline{Dipartimento di Matematica e Informatica}
		\centerline{Universit\`{a} di Cagliari}
		\centerline{Via Ospedale 72, 09124. Cagliari (Italy)}
		\medskip
	}
	\resetconstants{c}
%\tableofcontents
	\begin{abstract}
The following fully nonlinear attraction-repulsion and zero-flux chemotaxis model is studied: 
\begin{equation} \label{problem_abstract}
\tag{$\Diamond$}
\begin{cases}
u_t= \nabla \cdot \tonda*{(u+1)^{m_1-1}\nabla u -\chi u(u+1)^{m_2-1}\nabla v + \xi u(u+1)^{m_3-1}\nabla w} +\lambda u -\mu u^r  & \textrm{ in } \Omega \times (0,T_{max}),\\
\tau v_t=\Delta v-\phi(t,v)+f(u)  & \textrm{ in } \Omega \times (0,T_{max}),\\
\tau w_t= \Delta w - \psi(t,w) + g(u)& \textrm{ in } \Omega \times (0,T_{max}).
%u(x,0)=u_0(x), \; v(x,0)=v_0(x) & x \in \bar\Omega.
\end{cases}
%\showthe\font
\end{equation}
Herein, $\Omega$ is a bounded and smooth domain of $\R^n$, for $n\in \N$, $\chi,\xi, \lambda, \mu, r$ proper positive numbers, $m_1,m_2,m_3\in \R$, and $f(u)$ and $g(u)$ regular functions that generalise the prototypes $f(u) \simeq u^k$ and $g(u) \simeq u^l$, for some $k,l>0$ and all $u\geq 0$. Moreover, $\tau\in\{0,1\}$, and $\TM\in (0,\infty]$ is the maximal interval of existence of solutions to the model. Once suitable initial data $u_0(x),\tau v_0(x),\tau w_0(x)$ are fixed, we are interested in deriving sufficient conditions implying globality (i.e., $\TM=\infty$) and boundedness (i.e., $\norm{u(\cdot,t)}_{L^\infty(\Omega)}+\norm{v(\cdot,t)}_{L^\infty(\Omega)}+\norm{w(\cdot,t)}_{L^\infty(\Omega)}\leq C$ for all $t\in(0,\infty)$) of solutions to problem \eqref{problem_abstract}. This is achieved in these scenarios:
\begin{itemize}
\item [$\triangleright$] For $\phi(t,v)$ proportional to $v$ and $\psi(t,w)$ to $w$, whenever $\tau=0$ and provided one of the following conditions
\begin{enumerate}[label=(\Roman*),ref=\Roman*]
\item \label{El1} $m_2+k<m_3+l$,
\inlineitem \label{El2} $m_2+k<r$,
\inlineitem \label{El3} $m_2+k<m_1+\frac 2n$
\end{enumerate} 
is accomplished or $\tau=1$ in conjunction with one of these restrictions
\begin{enumerate}[label=(\roman*)]
\item  $\max{\graffa{m_2+k,m_3+l}}<r$, \hspace*{1cm}
\inlineitem $\max{\graffa{m_2+k,m_3+l}}<m_1+\frac 2n$,
\item $m_2+k<r$ and $m_3+l<m_1+\frac 2n$,
\inlineitem $m_2+k<m_1+\frac 2n$ and $m_3+l<r$;
\end{enumerate} 
\item [$\triangleright$] For $\phi(t,v)=\frac{1}{\abs*{\Omega}}\into f(u)$ and $\psi(t,w)=\frac{1}{\abs*{\Omega}}\into g(u)$, whenever $\tau=0$ if moreover one among \eqref{El1}, \eqref{El2}, \eqref{El3} is fulfilled.
\end{itemize}
Our research improves and extends some results derived in \cite{JiaoJadLi,GuoqiangBinATT-RepNonlinDiffSensLogistic,ChiyoYokotaBlow-UpAttRe, ColFraVig-ApplAnal}.
\end{abstract}
	
\section{Introduction and motivations}
\subsection{Some indications on attraction-repulsion chemotaxis models}\label{SectionAttr-Repulsion}

The general formulation of an attraction-repulsion chemotaxis model, incorporating logistic sources and linear or nonlinear productions, can be expressed as follows:
\begin{equation}\label{problemAttRep}
\hspace*{-0.2cm}
\begin{cases}
u_t= \nabla \cdot \left(D(u)\nabla u - S(u)\nabla v + T(u) \nabla w\right)+h(u)  & \text{ in } \Omega \times (0,\TM),\\
\tau v_t=\Delta v-\phi(t,v)+f(u)  & \text{ in } \Omega \times (0,T_{max}),\\
\tau w_t= \Delta w - \psi(t,w) + g(u)& \textrm{ in } \Omega \times (0,T_{max}),\\
%u_{\nu}=v_{\nu}=w_{\nu}=0 & \text{ on } \partial \Omega \times (0,\TM),\\
%u(x,0)=u_0(x), \; v(x,0)=v_0(x), \; w(x,0)=w_0(x)& x \in \bar\Omega,
%\end{cases}
%\\ & 
%\quad \textrm{with} \quad 
%\begin{cases}
u_{\nu}=v_{\nu}=w_{\nu}=0 & \text{ on } \partial \Omega \times (0,\TM),\\
u(x,0)=u_0(x),\; \tau v(x,0)=\tau v_0(x),\; \tau w(x,0)=\tau w_0(x)& x \in \bar\Omega,\\ %v(x,0)=v_0(x) & x \in \bar\Omega, \\  w(x,0)=w_0(x)& x \in \bar\Omega,
\end{cases}
\end{equation}
where $\Omega$ is a bounded and smooth domain in $\R^n$, $n\in \N$, and $D(u)$, $S(u)$, $T(u)$, $h(u)$, $\phi(t,v)$, $\psi(t,w)$, $f(u)$ and $g(u)$ are functions with specific regularity properties. Furthermore, $\tau \in\{0,1\}$ and additional regular initial data  $u_0(x)\geq 0$ and $ \tau v_0(x),\tau w_0(x) \geq 0$ are provided. The subscript $\nu$ in $(\cdot)_\nu$ denotes the outward normal derivative on $\partial \Omega$ and $\TM \in (0,\infty]$ represents the maximal temporal instant up to which solutions to the system do exist.

This model holds practical significance since it is applicable to studying inflammation in Alzheimer's disease. In this frame, microglia secrete both attractive and repulsive chemicals, and the system above describes the overall mechanisms of the involved quantities. In the specific linear diffusion, sensitivities and productions case (i.e., $D(u) = 1$, $S(u) = \chi u$,  $T(u) = \xi u $, and $f(u) = g(u) = u$), and in the absence of dampening external terms (i.e., $h(u)\equiv 0$), for $\phi(t,v) = v$ and $\psi(t,w) = w$
in \cite{Luca2003Alzheimer} the authors provide insights connected to gathering mechanisms for \eqref{problemAttRep} and develop numerical analyses within bounded intervals, particularly when $\tau = 0$.

\subsection{The attractive and the repulsive models}\label{SubsectionAttr-Repuls}
Model \eqref{problemAttRep} results from a combination of these perturbed signal-production mechanisms with aggregative effect   
\begin{equation}\label{problemOriginalKS} 
%\begin{cases}
u_t= \nabla \cdot(D(u)\nabla u- S(u)\nabla v)+ h(u)\quad \textrm{and} \quad 
\tau v_t=\Delta v- \phi(t,v)+f(u), \quad \textrm{ in } \Omega \times (0,\TM),
%\end{cases}
\end{equation}
and repulsive one 
\begin{equation}\label{problemOriginalKSCosnumption}
%\begin{cases}
u_t= \nabla \cdot(D(u)\nabla u+ T(u) \nabla w) + h(u) \quad \textrm{and} \quad 
\tau w_t=\Delta w-\psi(t,w)+g(u), \quad \textrm{ in } \Omega \times (0,\TM).
%\end{cases}
\end{equation} 
To properly understand model \eqref{problemAttRep}, it is essential to connect it with the biological mechanisms described in problems \eqref{problemOriginalKS} and \eqref{problemOriginalKSCosnumption}. The related partial differential equations are primarily employed to depict phenomena involving the spatial-temporal distribution of unicellular organisms (denoted as $u=u(x,t)$) within a confined and impenetrable environment ($\Omega$, with $(\cdot)_{\nu}=0$ on $\partial \Omega$). The movement of these organisms is influenced not only by natural diffusion (i.e., $\nabla \cdot D(u) \nabla u$), but also by the gradient of a chemical signal $v=v(x,t)$, referred to as a chemoattractant, which is produced at a rate $f(u)$ and leads to the aggregation of cells through the cross-diffusion $-\nabla \cdot S(u) \nabla v$. Additionally, a chemorepellent $w=w(x,t)$, secreted at a rate $g(u)$, contributes to the cell repulsion by means of the counterpart $+\nabla \cdot T(u) \nabla w$; moreover an external source $h(u)$ with both increasing and decreasing effects on the cell distribution takes part in the phenomena.
Naturally, in order to have well-posed systems, some initial configurations for cell density and chemical signals must be given; these are denoted by $u(x,0)=u_0(x)$, $\tau v(x,0)=\tau v_0(x)$, and $\tau w(x,0)= \tau w_0(x)$.

Regarding model \eqref{problemOriginalKS}, for the specific choice $D(u)=1$, $S(u)=\chi u$, $\phi(t,v)=v$, $f(u) = u$ and $h(u)=0$ (for which essentially the chemical signal $v$ increases with $u$), the natural homogenizing effect of the diffusion might be insufficient to make that the 
cell density equally disperses in the environment; indeed the drift term may force the system to experience a gathering process, resulting in the formation of highly concentrated spikes at some instants. This phenomenon, called \textit{chemotactic collapse} or \textit{blow-up at finite time}, is intimately related to the value $\chi m$ (being $m = \into u_0$ the initial mass of the particle distribution) and the dimension $n$. Mathematically, $\TM$ is finite and the solution $(u,v)$ becomes unbounded at $\TM$. Notably, when $n=1$, blow-up phenomena are excluded and in this case $\TM=\infty$ and the solution $(u,v)$ is bounded. However, for $n\geq 2$, chemotactic collapse occurs when $m \chi$ surpasses a critical value, herein denoted by $m_\chi$. If $m\chi$ is lower than $m_\chi$, no instability appears in the cells' motion. These findings are part of broader analyses, which explore existence and properties (globality, uniform boundedness, or blow-up) of solutions to the boundary-value problem associated to \eqref{problemOriginalKS}, especially in the parabolic-elliptic version ($\tau=0$). Further details can be found in references such as \cite{OsYagUnidim,HerreroVelazquez,JaLu,Nagai,WinklAggre} and others.

If in problem \eqref{problemOriginalKS} we fix $D(u)=1$, $S(u)=\chi u$, $\phi(t,v)=v$, $h(u)=0$, and $\tau=1$, by replacing the linear segregation $f(u) = u$ with a nonlinear one of the type $f(u) \simeq u^k$, with $0<k<\frac{2}{n}$ (for $n\geq 1$), implies (see \cite{LiuTaoFullyParNonlinearProd}) that all solutions remain bounded. On the contrary, when $\tau=0$ and $\phi(t,v)=\frac{1}{|\Omega|}\into u$, it is known that in spatially radial contexts (see \cite{WinklerNoNLinearanalysisSublinearProduction}), $\frac{2}{n}$ plays the role of critical value; indeed, solutions are bounded for any $n\geq 1$ and $0<k<\frac{2}{n}$, while blow-up phenomena may occur when $k>\frac{2}{n}$. (As it will be specified below, the term $\phi(t,v)=\frac{1}{|\Omega|}\into u$ influences on the production of $v$ in a nonlocal way: see Remark \ref{Rem:NonLoclDevi}.)

On the other hand, in the case of linear production rate, combining \eqref{problemOriginalKS} with terms representing population growth or decay, such as logistic sources (see \cite{verhulst}), seems quite natural. The equation for the particle density becomes 
\[u_t = \Delta u -\chi \nabla \cdot (u \nabla v) + h(u),\]
with $h(u)$ generally taking the form $h(u)=\lambda u- \mu u^r$, where $\lambda,\mu>0$ and $r>1$. The mathematical intuition suggests that the presence of a superlinear dampening effect may lead to boundedness and smoothness. This has been established for large $\mu$ (especially when $r=2$, $\phi(t,v)=v$, as shown in \cite{TelloWinkParEl}, \cite{W0} for respectively $\tau=0$ and $\tau=1$). However, for certain values of $r>1$, blow-up has been demonstrated in the parabolic-elliptic formulation, first with $\phi(t,v)=\frac{1}{|\Omega|}\into u$ for dimension $5$ or higher \cite{WinDespiteLogistic} (see also \cite{FuestCriticalNoDEA} for a recent improvement), and later even in three-dimensional domains, with $r<\frac{7}{6}$ and for the choice $\phi(t,v)=v$ (see \cite{Winkler_ZAMP-FiniteTimeLowDimension}). 

Finally, moving our attention to the situation where the produced signal $w$ has repulsive consequences on the cells'  motility, and this corresponds to model \eqref{problemOriginalKSCosnumption}, to the best of our knowledge no results regarding blow-up scenarios are available. This is likely due to the repulsive nature of the mechanism, and the existing literature on this topic is relatively limited, as seen in references like \cite{Mock74SIAM} and \cite{Mock75JMAA} which analyse similar contexts.   

\subsection{An overview on the state of the art for the linear version of model \eqref{problemAttRep}}\label{StateofTheArt}
What we have discussed for the purely attractive mechanism \eqref{problemOriginalKS} can be also observed for attractive-repulsive models with linear diffusion and drift terms. To be more precise, when in system \eqref{problemAttRep} we set $D(u)=1$, $S(u)=\chi u$, $T(u)=\xi u$, results dealing with blow-up and boundedness of solutions have been derived for specific expressions of the equations for $v$ and $w$. We mention two situations, known in the literature as local and nonlocal models.
\subsection*{The local model: $\phi(t,v)=\beta v$, $\psi(t,w)=\delta w$}
\vspace{2mm}
\begin{itemize}
\item [$\blacktriangleright$] \textit{The nonlogistic case ($h(u)=0$)}
\begin{itemize}
\item [$\triangleright$] The case $\tau=0$: In the case of linear growth for both the chemoattractant and the chemorepellent, i.e., $f(u) = \alpha u$, $\alpha>0$, and $g(u) = \gamma u$, $\gamma>0$, the difference $\Theta:=\chi \alpha-\xi\gamma$ between the parameters describing the impacts of the attraction and the repulsion, plays a crucial role. Specifically, when $\Theta<0$ (indicating a regime where the repulsion dominates the attraction), all solutions to the model are globally bounded in any dimension. Conversely, when $\Theta>0$ (emphasizing now, that the attraction is the dominating effect) and $n=2$, unbounded solutions can be detected (for further details, see \cite{GuoJiangZhengAttr-Rep,LI-LiAttrRepuls,TaoWanM3ASAttrRep,VIGLIALORO-JMAA-BlowUp-Attr-Rep,YUGUOZHENG-Attr-Repul}). On the other hand, for more general expressions of the production laws $f$ and $g$, and more exactly those generalising the prototypes $f(u) \simeq \alpha u^k$, $k>0$, and $g(u)\simeq \gamma u^l$, $l>0$, the following recent results valid for $n\geq 2$ apply (\cite{ViglialoroMatNacAttr-Repul}):
For any $\alpha,\beta,\gamma,\delta,\chi>0$, and $l>k\geq 1$ (or $k>l\geq 1$), there exists $\xi^*>0$ (or $\xi_*>0$) such that if $\xi>\xi^*$ (or $\xi\geq \xi_*$), any sufficiently regular initial datum $u_0(x)\geq 0$ (or $u_0(x)\geq 0$ small in some Lebesgue space) leads to a unique classical solution which remains uniformly bounded in time.    
Moreover, the same conclusion holds for any $\alpha, \beta, \gamma, \delta, \chi, \xi>0$, and any sufficiently regular $u_0(x)\geq 0$ if the conditions $0<k<1$ and $l=1$ are satisfied (see also \cite{ColFraVig-ApplAnal} for further improvements in some situations).
 
\item [$\triangleright$] The case $\tau=1$: The research in \cite{TaoWanM3ASAttrRep} establishes that in two-dimensional domains, when $f(u) = \alpha u$ and $g(u) = \gamma u$, sufficiently smooth initial data lead to solutions that are globally bounded in time. This holds true if the condition $\Theta<0$ is met (where $\Theta=\chi \alpha-\xi \gamma$) and either  
\[
\beta=\delta \quad \text{or} \quad  -\frac{\chi^2\alpha^2(\beta-\delta)^2}{2\Theta \beta^2} \int_\Omega u_0(x)\leq C(\Omega), \quad \textrm{for some constant }\; C(\Omega)>0.
\]
Additionally, in three-dimensional ball-shaped domains, \cite{LankeitJMAA-BlowupParabolicoAttr-Rep} highlights the occurrence of blow-up at a finite time under the condition $\Theta = \chi \alpha - \xi \gamma>0$. 
\end{itemize}  
\vspace{2mm}
\item [$\blacktriangleright$] \textit{The logistic case ($h(u)=\lambda u -\mu u^r$, $\lambda,\mu>0$, $r>1$)}
\begin{itemize}
\item [$\triangleright$] {The case $\tau=0$}: For linear rates of the chemoattractant and the chemorepellent, $f(u) = \alpha u$, $\alpha>0$, and $g(u) = \gamma u$, $\gamma>0$, finite time blow-up is proved in \cite{ChiyoMarrasTanakaYokota2021}
under the assumption $\Theta=\chi \alpha -\xi \gamma >0$ in $n$-dimensional balls ($n \geq 3$), for some $r=r(n)$ close to $1$ (additionally, also estimates of the blow-up time are given). Otherwise, for nonlinear behaviour of $f$ and $g$ ($f(u) \simeq u^k$, $g(u)\simeq u^l$), in \cite{LiangEtAlAtt-RepNonLinProdLogist-2020} the authors show, inter alia, that if $k<\max\{l, r-1,\frac{2}{n}\}$ all solutions are globally bounded; the related long-time behaviour of these solutions is studied in  \cite{XinluEtAl2022-Asymp-AttRepNonlinProd}.
\item [$\triangleright$] {The case $\tau=1$}: For $r=2$ and in the three-dimensional setting whenever $\mu$ is sufficiently large, boundedness and rate of convergence to constant equilibria are discussed in \cite{GuoqiangBin-2022-3DAttRep}.
\end{itemize}
\end{itemize}
\subsection*{The nonlocal model: $\tau=0$, $\phi(t,v)=\frac{1}{\abs{\Omega}}\into f(u)$, $\psi(t,w)=\frac{1}{\abs{\Omega}}\into g(u)$}
In \cite{LiuLiBlowUpAttr-Rep} it is proved, together with other results, that 
for $f(u) \simeq \alpha u^k$ and $g(u)\simeq \gamma u^l$ with $k>\frac{2}{n}$ and $k>l$, unbounded solutions can be detected, in absence of logistics. 
By contrast, also in the case $h\equiv 0$, in \cite{ColFraVig-ApplAnal} the authors show the boundedness of solutions provided that $k<l$ or $k=l$ and $\Theta:=\chi \alpha - \xi \gamma <0$, or $l=k \in (0,\frac{2}{n})$ and $\Theta\geq0$.
\vspace{0.5cm}

So far, we only restricted our discussion to linear diffusion and sensitivities in model \eqref{problemAttRep}. Conversely, this project aims at providing results for alike models (both including local and nonlocal effects) involving nonlinear diffusion and drift terms, whose prototypes are $D(u)\simeq u^{m_1}, S(u)\simeq u^{m_2}, T(u)\simeq u^{m_3}$, where $m_1,m_2,m_3\in\R$. 
\section{The nonlinear local and nonlocal models, known results and main claims}
%\subsection{Definition of the models}
From the pure mathematical point of view, our attention is direct to these two systems, one of local type (indicated with ($\mathcal{L}$)) and the other of nonlocal nature, ($\mathcal{N}$):
	\begin{equation}\label{eq:localproblem}
		\begin{dcases}\tag{$\mathcal{L}$}
			u_t= \nabla \cdot \tonda*{(u+1)^{m_1-1}\nabla u -\chi u(u+1)^{m_2-1}\nabla v + \xi u(u+1)^{m_3-1}\nabla w} +\lambda u -\mu u^r & \text{in } \Omega \times (0,\TM),\\
			\tau v_t= \Delta v - \beta v  +f(u) & \text{in } \Omega \times (0,\TM),\\
			\tau w_t= \Delta w -  \delta w + g(u) & \text{in } \Omega \times (0,\TM),\\
			u_{\nu}=v_{\nu}=w_{\nu}=0 & \text{on } \partial \Omega \times (0,\TM),\\
			u(x,0)=u_0(x), \tau v(x,0)=\tau v_0(x), \tau w(x,0)=\tau w_0(x)  & x \in \bar\Omega,
		\end{dcases}
	\end{equation}
	and
	\begin{equation}\label{eq:nonlocalproblem}
		\begin{dcases}\tag{$\mathcal{N}$}
			u_t= \nabla \cdot \tonda*{(u+1)^{m_1-1}\nabla u -\chi u(u+1)^{m_2-1}\nabla v + \xi u(u+1)^{m_3-1}\nabla w}+\lambda u -\mu u^r & \text{in } \Omega \times (0,\TM),\\
			0= \Delta v - \frac{1}{\abs{\Omega}}\into f(u)  + f(u) & \text{in } \Omega \times (0,\TM),\\
			0= \Delta w - \frac{1}{\abs{\Omega}}\into g(u)  +g(u) & \text{in } \Omega \times (0,\TM),\\
			u_{\nu}=v_{\nu}=w_{\nu}=0 & \text{on } \partial \Omega \times (0,\TM),\\
			u(x,0)=u_0(x) & x \in \bar\Omega,\\
			\int_\Omega v(x,t)\,dx=\int_\Omega w(x,t)\,dx=0 & \text{for all } t\in  (0,\TM).
		\end{dcases}
	\end{equation}
 \begin{remark}\label{Rem:NonLoclDevi}
 In problem \eqref{eq:localproblem}, the equation $\tau v_t=\Delta v-\beta v+f(u)$ describes the local interaction between the involved quantities at each spatial point $x$; contrarily,  in \eqref{eq:nonlocalproblem} the nonlocal term $\frac{1}{|\Omega|}\int_\Omega f(u)$ appears, and it takes into account the entire distribution of $u$ all over the domain $\Omega$. (Naturally, the same observation can be done for the equations related to $w$.) 
 
Additionally,  let us clarify that in the nonlocal model,  $v$ and $w$ stand for the deviations of the chemoattractant and the chemorepellent (and not for the chemicals themselves). Since by definition the deviation is the difference between the signal concentration and its mean value,  it follows that the means of $v$ and $w$ vanish, exactly as  specified in the last positions of problem \eqref{eq:nonlocalproblem}. As a consequence, $v$ and $w$ change sign. 
 \end{remark}
Consistently with the above problems, this paper is contextualized in the frame of a series of results dealing with local and nonlocal, and linear and nonlinear, attraction-repulsion chemotaxis systems; in particular this research improves and extends already known analyses in the literature (we will give more details in the sequel).

In the specific, as far as we know, for such discussed fully nonlinear versions, the literature is rather poor; nevertheless, issues dealing with boundedness and blow-up of solutions have been addressed. 
More precisely, for problem \eqref{eq:localproblem} we mention the works 
\cite{GuoqiangBinATT-RepNonlinDiffSensLogistic} and \cite{ChiyoYokotaBlow-UpAttRe}, where in the first the authors ensure the boundedness of solutions, for both $\tau=0$ and $\tau=1$, under suitable assumptions on the data; in the second boundedness and blow-up analyses, in this case for the simplified parabolic-elliptic-elliptic version, are discussed. For the nonlocal problem \eqref{eq:nonlocalproblem}, and in the context of nonlinear productions, in \cite{wang2023blow} blow-up scenarios have been detected  for $m_1\in\R$ and $m_2=m_3=1$; in turn, in the recent paper \cite{ColFraVig-StudApplMat} this result has been extended to the case $m_1\in\R$ and $m_2=m_3\in \R$. As to the boundedness counterpart, at this point we simply make reference to \cite{JiaoJadLi}.   
(Since in this research we will improve some results derived in \cite{JiaoJadLi,GuoqiangBinATT-RepNonlinDiffSensLogistic,ChiyoYokotaBlow-UpAttRe, ColFraVig-ApplAnal},  we will spend more words on these papers in the below Remark \ref{RemNonLoc} and Remark \ref{RemLoc}.)
\subsection{Presentation of the main theorems}
In order to present our claims, some positions have to be previously fixed. In particular, here we assume that $\Omega$, sources $f, g$ and the initial data $u_0=u_0(x), \tau v_0=\tau v_0(x)$ and $\tau w_0=\tau w_0(x)$  are such that 
\begin{equation}\label{eq:reglocal}
	\begin{cases}
            \Omega\subset\R^n, n\in\N, \textrm{ is a bounded domain with smooth boundary $\partial \Omega$},\\
			f,g: [0,\infty) \rightarrow \R^+, \textrm{ with } f,g\in C^1([0,\infty)), 
			\\ 
			u_0, \tau v_0, \tau w_0: \bar{\Omega}  \rightarrow \R^+, \textrm{ with } u_0,  \tau v_0,   \tau w_0 \in W^{1,\infty}(\Omega); 
	\end{cases}
	\end{equation}
	moreover, we suppose, for $\alpha,l,k>0$ and $0<\gamma_0\leq\gamma_1$,
	\begin{equation}\label{eq:disf}
		0\leq f(s)\leq \alpha (s+1)^k \quad \text{and} \quad  \gamma_0(s+1)^l\leq g(s)\leq \gamma_1 (s+1)^l.
	\end{equation}
 Additionally, we will frequently invoke these
	\begin{assumptions}\label{generalassumptions}
		For $n\in\N$, $m_1,m_2,m_3\in\R$, $k,l>0$ and $r>1$ let us set:
		\begin{enumerate}[label={($\mathcal{A}_{\arabic*}$)}]
			\item \label{M1}  $m_2+k<m_3+l$;
			\inlineitem \label{M2}  $m_2+k<r$;
			\inlineitem \label{M3}  $m_2+k<m_1+\frac 2n$;
			\inlineitem \label{M4}  $m_3+l<r$;
			\inlineitem \label{M5}  $m_3+l<m_1+\frac 2n$.
		\end{enumerate}
	\end{assumptions}
\noindent The above preparations allow us to give the following two theorems.
	\begin{theorem}\label{theo:localtau0}
		For $\tau=0$, let $\Omega$, $f,g,u_0$ comply with hypotheses in \eqref{eq:reglocal}, \eqref{eq:disf}, $\lambda,\mu,\chi,\xi,\beta,\delta>0$ and $q>n$.  Additionally, let one among \ref{M1}, \ref{M2}, \ref{M3} in Assumptions \ref{generalassumptions} hold true. Then problems \eqref{eq:localproblem} and \eqref{eq:nonlocalproblem} admit a unique solution
        \[
        u \in C^0(\bar{\Omega} \times [0,\infty)) \cap C^{2,1}(\bar{\Omega} \times (0,\infty)) \;\textrm{and}\;
        v, w \in C^{2,0}(\bar{\Omega} \times (0,\infty)) \cap
        L^{\infty}_{loc}((0,\infty); W^{1,q}(\Omega)), 
        \]
		such that $u$ is nonnegative and $u,v,w$ are bounded on $\bar\Omega \times [0,\infty)$. In particular, $v,w$ are as well nonnegative on $\bar\Omega \times [0,\infty)$ for problem \eqref{eq:localproblem}.
	\end{theorem}
 
	\begin{theorem}\label{theo:localtau1}
		For $\tau=1$, let the remaining hypotheses of Theorem \ref{theo:localtau0} be satisfied.  Additionally, let one among \ref{M2} and \ref{M3} jointly with one between \ref{M4} and \ref{M5} in Assumptions \ref{generalassumptions} hold true. Then problem \eqref{eq:localproblem} admits a unique solution 
		\begin{equation*}
		  u \in C^0(\bar{\Omega} \times [0,\infty)) \cap C^{2,1}(\bar{\Omega} \times (0,\infty)) \;\textrm{and}\;
		v, w \in C^0(\bar{\Omega} \times [0,\infty)) \cap C^{2,1}(\bar{\Omega} \times (0,\infty)) \cap L^{\infty}_{loc}((0,\infty); W^{1,q}(\Omega)), 
		\end{equation*}
		such that $u,v,w$ are nonnegative and bounded on $\bar\Omega \times [0,\infty)$.
	\end{theorem}
 As anticipated in the introductory part, let us specify in which sense our results improve what known so far in the literature.
 \begin{remark}\label{RemNonLoc}{\normalfont(Comparisons with nonlocal problem \eqref{eq:nonlocalproblem})}
 Let us make these observations.
\begin{itemize}
    \item[$\blacktriangleright$] Nonlinear diffusion and sensitivities ($m_1, m_2, m_3 \in \R$): in \cite[Theorem 1.1]{JiaoJadLi} boundedness for problem \eqref{eq:nonlocalproblem} is achieved for each of the following cases:
    \begin{enumerate}[label=\roman*)] 
    \item \label{Li1} $m_3+l<m_2+k+1$ and $m_2+k<m_1+\frac{2}{n}$;
    \inlineitem \label{Li2} $m_2+k<m_3+l<m_1+\frac{2}{n}$;
    \inlineitem \label{Li3} $\max\{m_2+k, m_3+l\}<r<m_1+\frac{2}{n}$.
    \end{enumerate}
    It is seen that assumption \ref{M1} in Theorem \ref{theo:localtau0} is sharper with respect to \ref{Li2}, and \ref{M2} with respect to \ref{Li3} and, finally,  
    \ref{M3} with respect to \ref{Li1}.
    \item[$\blacktriangleright$] Linear diffusion and sensitivities ($m_1=m_2=m_3=1$): in \cite[Theorem 2.3]{ColFraVig-ApplAnal} boundedness for problem \eqref{eq:nonlocalproblem} is established, among other situations, whenever:
    \begin{enumerate}[label=(\alph*)]
\item \label{CFV1} $k<l$;  
\inlineitem \label{CFV3} $k=l \in (0,\frac{2}{n})$ and $\Theta_0:=\chi \alpha - \xi \gamma_0\geq 0$. 
\end{enumerate} 
Evidently, if from the one hand \ref{CFV1} is recovered from \ref{M1}, on the other hand, also in this case some milder conclusions can be observed: in the specific \ref{M3} improves \ref{CFV3} and provides a further situation where $l>0$ and $k\in \left(0,\frac{2}{n}\right)$.
\end{itemize}
\end{remark}
\begin{remark}\label{RemLoc}{\normalfont(Comparisons with local problem \eqref{eq:localproblem})}
Even for the local problem, the present analysis gives some more precise insight in the context of attraction-repulsion chemotaxis mechanisms.   
\begin{itemize}
    \item[$\blacktriangleright$] Nonlinear diffusion and sensitivities ($m_1, m_2, m_3 \in \R$): 
\begin{itemize}
    \item[$\triangleright$] $\tau=0$: In \cite[Theorem 3.1 and Theorem 3.5]{ChiyoYokotaBlow-UpAttRe}, where the authors consider the problem with linear productions ($k=l=1$), globality of solutions is carried out for the case $m_2 \leq m_3$, so leaving room to the case $m_2>m_3$; in particular, for either $m_2 < m_3$ or $m_2=m_3$ and $\Theta_0=\chi \alpha - \xi \gamma_0<0$ (naturally, for $k=l=1$ we have $\gamma_0=\gamma_1$), boundedness is, respectively, achieved. Indeed, assumptions \ref{M2} and \ref{M3} of Theorem \ref{theo:localtau0} cover the case $m_2>m_3$.

Furthermore, Theorem \ref{theo:localtau0} extends \cite[Theorem 1.3]{GuoqiangBinATT-RepNonlinDiffSensLogistic}, where only the case $m_1 \in \R$, $m_2=m_3=1$ and $k,l \geq 1$ is addressed; particularly assumption \ref{M3} is less restricted than \cite[Theorem 1.3, (i)]{GuoqiangBinATT-RepNonlinDiffSensLogistic}.
\item[$\triangleright$] $\tau=1$: In \cite[Theorem 1.1]{GuoqiangBinATT-RepNonlinDiffSensLogistic} boundedness of solutions is established, inter alia, if   
\begin{equation}\label{RL}
\max\{m_2+k, m_3+l\}<m_1+\frac{2}{n} \quad \text{or} \quad \max\{m_2+k, m_3+l\}<r
\end{equation}
are fulfilled. In this way, (\ref{M2}, \ref{M4}) and (\ref{M3}, \ref{M5}) of Theorem \ref{theo:localtau1} coincide with \eqref{RL}, but (\ref{M2}, \ref{M5}) and (\ref{M3}, \ref{M4}) provide a further scenario towards boundedness.
\end{itemize}
    \item[$\blacktriangleright$] Linear diffusion and sensitivities ($m_1=m_2=m_3=1$):
 For $\tau=0$ if we compare Theorem \ref{theo:localtau0} with \cite[Theorem 2.1]{ColFraVig-ApplAnal}, not only the same conclusion given in the second item of Remark \ref{RemNonLoc} still applies but even \cite[Thorem 2.1, (2)]{ColFraVig-ApplAnal} is turned into a weaker condition. For the fully parabolic case, assumption \ref{M3} and \ref{M5} of Theorem \ref{theo:localtau1} imply $l,k \in (0,\frac{2}{n})$, so improving \cite[Theorem 2.2]{ColFraVig-ApplAnal}.  
\end{itemize}
\end{remark}
\section{Adjusting parameters and recalling necessary results}
Let us now dedicate to summarise some tools which will be used in our reasoning. These are connected to algebraic inequality and regularity results for Partial Differential Equations. 
\begin{lemma}\label{lem:disab}
Let $A,B \geq 0$ and $d_1, d_2>0$ be such that $d_1 + d_2 <1$. Then for all $\epsilon >0$ there exists $c=c_{\varepsilon}>0$ such that
	\[
	A^{d_1}B^{d_2} \leq \epsilon(A+B) +c. 
	\]
\begin{proof}
The proof can be found in \cite[Lemma 4.3]{frassuviglialoro1}.
\end{proof}
\end{lemma}
	\begin{lemma} \label{eq:regularity}
		Let $\Omega\subset\R^n$ satisfy condition in \eqref{eq:reglocal}, $\eta>0$ and $q>\max\{n,\frac{1}{\eta}\}$.
		\begin{itemize}
			\item [$\diamond$] {\bf\emph{Elliptic regularity}:} Let $\vartheta \in (0,1)$. If $\psi\in C^{\vartheta}(\bar\Omega)$, then the solution $z\in C^{2+\vartheta}(\bar\Omega)$  of
			\begin{equation*}
				\begin{cases}
					-\Delta z + \eta z=\psi  & \text{in } \Omega,\\
					z_{\nu}=0 & \text{on } \partial\Omega,
				\end{cases}
			\end{equation*}
			is such that 
			\begin{equation}\label{eq:tau0}
            \nabla z \in L^\infty(\Omega).
            \end{equation}
   Moreover for any $\hat{c},\rho>0$, there is $C_\rho>0$ such that
\begin{equation}\label{ellipticreg}
    \hat{c}\into z^q\leq \rho \into \psi^q + C_\rho\tonda*{\into \psi}^q.
\end{equation}
			\item [$\diamond$] {\bf\emph{Parabolic regularity}:} Let $T\in(0,\infty]$. Then, for $\psi \in L^q([0,T);L^q(\Omega))$ and $z_0\in W^{2,q}(\Omega)$ with $\partial_{\nu} z_0 = 0$ on $\partial\Omega$, there exists $C_P=C_P(\Omega, q, \rVert z_0\rVert_{W^{2,q}(\Omega)})$ such that every solution $z\in W^{1,q}_{loc}([0,T);L^q(\Omega))\cap L^q_{loc}([0,T);W^{2,q}(\Omega))$ of 
			\begin{equation*}
				\begin{cases}
					z_t= \Delta z - \eta z  +\psi & \text{in } \Omega \times (0,T),\\
					\partial_{\nu} z = 0 & \text{on } \partial\Omega \times (0,T),\\
					z(\cdot,0)=z_0 & \text{on } \Omega,
				\end{cases}
			\end{equation*}
			satisfies
	\begin{equation}\label{eq:tau1}
				\int_0^t e^s \into \abs*{\Delta z(\cdot,s)}^q\,ds \leq C_{P}\quadra*{1+\int_0^t e^s \into\abs{\psi(\cdot,s)}^q\,ds} \quad \text{for all } t\in(0,T).
			\end{equation}
			Additionally, if  $\psi \in L^\infty([0,T);L^q(\Omega))$ then
			\begin{equation}\label{eq:tau1extension}
				\nabla z \in L^\infty((0,T);L^\infty(\Omega)).
			\end{equation}
		\end{itemize}
		\begin{proof} Relation \eqref{eq:tau0} is an obvious consequence of \cite[Theorem IX.33]{BrezisBook}. For estimate \eqref{ellipticreg}, we indicate \cite[Lemma 2.2]{WinklerHowFar} (and also \cite[Lemma 3.1]{ViglialoroMatNacAttr-Repul}) and we precise that the power $q$ may be replaced by any other one larger than $1$. (Since $q$ will be chosen arbitrarily large we preferred to synthesize the nomenclature.) On the other hand, for \eqref{eq:tau1} we refer to \cite[Lemma 3.6]{IshidaLankeitVigliloro-Gradient}, whereas for \eqref{eq:tau1extension} we can invoke \cite[Lemma 4.1]{HorstWink} and the embedding $W^{1,\frac{n q}{(n-q)_+}}(\Omega)\subset L^\infty(\Omega)$, valid for $q>n$. (The reason why we need $q>\frac{1}{\eta}$ is understandable seeing  \cite[Lemma 3.6]{IshidaLankeitVigliloro-Gradient}.)
		\end{proof}
	\end{lemma}

The following lemma defines crucial parameters; its proof is based on some of the relations fixed in Assumptions \ref{generalassumptions}.
	\begin{lemma}\label{lem:p}
		Let $n,k,l,m_1,m_2,m_3$ be as in Assumptions \ref{generalassumptions} and let relations \ref{M3} and \ref{M5} be valid. Then, there exists $\Bar{p}>1$ such that, for all $p>\Bar{p}$, $q>1$ and
		% all $p>\max{\{1,\frac n2,2-(m_2+\alpha), \frac n\gamma (m_2+\alpha-\gamma)-(m_2+\alpha-1)\}}$ 
		\begin{equation*}
			\theta(p)\coloneqq\frac{\frac{p+m_1-1}{2}-\frac{p+m_1-1}{2(p+m_2+k-1)}}{\frac{p+m_1-1}{2}-\frac 12+\frac 1n}, \quad %\quad 
			\sigma(p)\coloneqq\frac{2(p+m_2+k-1)}{p+m_1-1}, \quad \theta_1(p)\coloneqq \frac{\frac{p+m_1-1}{2}-\frac{p+m_1-1}{2(p+m_3+l-1)}}{\frac{p+m_1-1}{2}-\frac 12+\frac 1n}, \quad %\quad 
			\sigma_1(p)\coloneqq\frac{2(p+m_3+l-1)}{p+m_1-1}
		\end{equation*}
		\begin{equation*} 
			\theta_2(p)\coloneqq \frac{\frac{p+m_1-1}{2}-\frac{p+m_1-1}{2l}}{\frac{p+m_1-1}{2}-\frac 12+\frac 1n}, \quad \quad 
			\theta_3(p)\coloneqq \frac{\frac{p+m_1-1}{2}-\frac{p+m_1-1}{2p}}{\frac{p+m_1-1}{2}-\frac 12+\frac 1n}, \quad \quad
			\theta_4(p)\coloneqq\frac{\frac{p+m_1-1}{2}-\frac{p+m_1-1}{2q}}{\frac{p+m_1-1}{2}-\frac 12+\frac 1n}, \quad \quad 
			\sigma_2(p)\coloneqq\frac{2(p+q)}{p+m_1-1},
		\end{equation*}
		these relations hold
		\begin{table}[H]
			\centering
			\begin{subequations}
				\begin{subtable}[h]{0.24\textwidth}
					\centering
					\begin{equation}\label{eq:theta}
						0<\theta<1,
					\end{equation}
				\end{subtable}
				\hfill
				\begin{subtable}[h]{0.24\textwidth}
					\centering
					\begin{equation}\label{eq:sigmatheta}
						0<\frac{\sigma\theta}{2}<1,
					\end{equation}
				\end{subtable}
				\hfill
				\begin{subtable}[h]{0.24\textwidth}
					\centering
					\begin{equation}\label{eq:thetabar}
						0<\theta_1<1,
					\end{equation}
				\end{subtable}
				\hfill
				\begin{subtable}[h]{0.24\textwidth}
					\centering
					\begin{equation}\label{eq:sigmathetabar}
						0<\frac{\sigma_1\theta_1}{2}<1,
					\end{equation}
				\end{subtable}
				\\
				\begin{subtable}[h]{0.18\textwidth}
					\centering
					\begin{equation}\label{eq:thetahat}
						0<\theta_2<1,
					\end{equation}
				\end{subtable}
				\hfill
				\begin{subtable}[h]{0.20\textwidth}
					\centering
					\begin{equation}\label{eq:sigmathetahat}
						0<\frac{\sigma_1\theta_2}{2}<1,
					\end{equation}
				\end{subtable}
				\hfill
				\begin{subtable}[h]{0.18\textwidth}
					\centering
					\begin{equation}\label{eq:thetatilde}
						0<\theta_4<1,
					\end{equation}
				\end{subtable}
				\hfill
				\begin{subtable}[h]{0.20\textwidth}
					\centering
					\begin{equation}\label{eq:sigmathetatilde}
						0<\frac{\sigma_2\theta_4}{2}<1,
					\end{equation}
				\end{subtable}
				\hfill
				\begin{subtable}[h]{0.18\textwidth}
					\centering
					\begin{equation}\label{eq:thetaunder}
						0<\theta_3<1,
					\end{equation}
				\end{subtable}
			\end{subequations}
		\end{table}
\noindent where \eqref{eq:thetahat} and \eqref{eq:sigmathetahat} are valid only for $l>1$.
  \begin{proof}
      First let us consider the function $\theta(p)$. Since
      \begin{math}
          \lim_{p\to\infty} \theta(p) = 1,
      \end{math}
     we have that $\theta(p)$ is definitively positive; on the other hand, $\theta(p)$ increases for $p$ sufficiently large so that 1 is an upper bound. Therefore, \eqref{eq:theta} is proved. Relations \eqref{eq:thetabar}, \eqref{eq:thetaunder}, \eqref{eq:sigmatheta} and \eqref{eq:sigmathetabar} follow by means of analogous arguments, once for \eqref{eq:sigmatheta} and \eqref{eq:sigmathetabar} conditions \ref{M3} and \ref{M5} are respectively taken into account.
      Let us consider now the functions $\theta_2(p)$ and $\sigma_1(p)$. Relations \eqref{eq:thetahat} and \eqref{eq:sigmathetahat} are consequence of
      \begin{equation*}
          \lim_{p\to\infty} \theta_2(p) = \lim_{p\to\infty} \frac{\sigma_1(p)\theta_2(p)}{2}= 1 - \frac1l \in (0,1),  \quad \text{for $l > 1$}.
      \end{equation*}
      In the same way, relations \eqref{eq:thetatilde} and \eqref{eq:sigmathetatilde} are given for every $q>1$.

      From all of the above, it is possible to find $\Bar{p}>1$ such that the previous relations are complied for every $p>\Bar{p}$.
  \end{proof}
	\end{lemma}
	
\section{Local existence and extensibility criterion}
A first necessary step on which our computations have to rely is the local-in-time existence of classical solutions to systems \eqref{eq:localproblem} and \eqref{eq:nonlocalproblem}. The succeeding requirement is, indeed, providing a criterion capable to turn such solutions into global ones. The following lemma focuses on these two aspects.
\begin{lemma}\label{lem:existence}
For $\tau \in \{0,1\}$, let $\Omega$, $f,g,u_0,\tau v_0$ and $\tau w_0$ comply with hypotheses in \eqref{eq:reglocal}. Moreover, let $\chi,\xi,\beta,\delta, \lambda,\mu>0$, $m_1,m_2,m_3\in \R$, $r>1$, and $q>n$. Then there exist $\TM \in (0,\infty]$ and a unique solution $(u,v,w)$ to problems \eqref{eq:localproblem} and \eqref{eq:nonlocalproblem}, defined in $\Omega \times (0,\TM)$ and such that if $\tau=0$
\begin{equation*}
u \in C^0(\bar{\Omega} \times [0,\TM)) \cap C^{2,1}(\bar{\Omega} \times (0,\TM)) \;\textrm{and}\;
v, w \in C^{2,0}(\bar{\Omega} \times (0,\TM)) \cap L^{\infty}_{loc}((0,\TM); W^{1,q}(\Omega)),
\end{equation*}
whereas if $\tau = 1$
\begin{equation*}
	 u \in C^0(\bar{\Omega} \times [0,\TM)) \cap C^{2,1}(\bar{\Omega} \times (0,\TM)) \;\textrm{and}\;
	v, w \in C^0(\bar{\Omega} \times [0,\TM)) \cap C^{2,1}(\bar{\Omega} \times (0,\TM)) \cap L^{\infty}_{loc}((0,\TM); W^{1,q}(\Omega)). 
\end{equation*}
The components $(u,v,w)$ of solutions to problem \eqref{eq:localproblem} are nonnegative, whereas for \eqref{eq:nonlocalproblem} only $u$ is nonnegative. 
 
In addition,
\begin{equation}\label{eq:extensibility}
\text{if} \quad \TM<\infty \quad \text{then} \quad \limsup_{t \to \TM} \|u(\cdot,t)\|_{L^{\infty}(\Omega)}=\infty.
\end{equation}
Finally, let $u\in L^\infty((0,\TM),L^p(\Omega))$ for all $p>1$; then $u\in L^\infty((0,\TM),L^\infty(\Omega))$ and $\TM=\infty$.
\begin{proof}
Existence and uniqueness of solutions can be established using well known techniques based on fixed point arguments and elliptic and parabolic regularity results: we refer to \cite{cieslak2008finite, nagai1995blow,
wang2016quasilinear, winkler2010boundedness}. 

Let us spend some words on the last implication. Naturally $u = u(x,t)$, with $(x,t)\in \Omega\times(0,\TM)$, also classically solves problem (A.1) in \cite[Appendix A]{TaoWinkParaPara} for 
\begin{equation*}
D(x,t,u)=(u+1)^{m_1-1},\quad \Tilde{f}(x,t)=-\chi u(u+1)^{m_2-1}\nabla v+\xi  u(u+1)^{m_3-1}\nabla w,\quad g(x,t)=\lambda u-\mu u^r. 
\end{equation*}
Specifically, by making use of the Neumann boundary conditions, we can see that (A.2)--(A.5) are complied. On the other hand, for any $\lambda,\mu>0$ and $r>1$, it holds that $\lambda u-\mu u^{r}$ has a positive maximum $L$ at $u_M=\left(\frac{\lambda}{r \mu}\right)^{\frac{1}{r-1}}$, so that from $g(x,t)\leq L$ in $\Omega \times (0,T_{max})$ the second inclusion (A.6) is obtained for any $q_2>1$. 
Since, by hypotheses, $u\in L^\infty((0,\TM);L^p(\Omega))$ for every $p > 1$, we have that $\into u^p$ is uniformly bounded on $(0,\TM)$ for $p$ arbitrary large (without relabelling it) and henceforth conditions (A.7)--(A.10) are fulfilled. As to the first inclusion of \cite[(A.6)]{TaoWinkParaPara}, if $\tau=0$ in problem \eqref{eq:localproblem} (or model \eqref{eq:nonlocalproblem}), from the gained inclusion $u\in L^\infty((0,\TM);L^p(\Omega))$ we have $f(u)\in L^\infty((0,\TM);L^p(\Omega))$, and in turn relation \eqref{eq:tau0} provides  $\nabla v\in L^\infty((0,\TM);L^\infty(\Omega))$, and similarly $\nabla w \in L^\infty((0,\TM);L^\infty(\Omega))$. For $\tau = 1$ we directly invoke \eqref{eq:tau1extension}; in both cases we have that, for any $q_1>1$,
$\Tilde{f} \in L^\infty((0,\TM);L^{q_1}(\Omega))$. As a consequence of what explained, we have the claim by virtue of \cite[Lemma A.1]{TaoWinkParaPara} and the extensibilty criterion \eqref{eq:extensibility}.
\end{proof}
 \end{lemma}
\section{Some a priori estimates}
In this section we will dedicate to derive some uniform-in-time estimates of the previously  gained local solution the investigated problems. In this frame, from now on we will tacitly assume that
\begin{itemize}
    \item  [$\triangleright$] all the constants $c_i$ ($i=1,2,\ldots$) appearing below are positive,
    \item [$\triangleright$] the triplet $(u,v,w)$ indicates the local solution to models \eqref{eq:localproblem} or \eqref{eq:nonlocalproblem} (naturally recognizable from the context), obtained in Lemma \ref{lem:existence}.
\end{itemize}
The forthcoming lemmas hold for models \eqref{eq:localproblem} and \eqref{eq:nonlocalproblem}.
\begin{lemma}\label{boundednessMass}
The component $u$ is such that the associated mass $\int_\Omega u(x,t)dx$ is uniformly bounded over $(0,\TM)$; more specifically, 
\begin{equation*}
  \int_{\Omega} u \leq M:=\max\left\{\displaystyle\int_\Omega u_0(x)\,dx, \left(\frac{\lambda}{\mu} |\Omega|^{r-1}\right)^{\frac{1}{r-1}} \right\} \quad \text{for all } t \in [0,\TM). 
\end{equation*} 
\begin{proof}
This property can be proved by integrating over $\Omega$ the first equation of the models and then by applying the H\"{o}lder inequality and ODI comparison principles.
\end{proof}
\end{lemma}
\begin{lemma}\label{lem:GN+YOUNG}
Let hypotheses of Lemma \ref{lem:p} be valid and $\Bar{p}>1$ the value therein found. Then for every $\hat{c},\varepsilon_1,\varepsilon_2>0$ and for every $p>\Bar{p}$, there exist some constants $\const{a2},\const{a4},\const{a6},\const{a4s},\const{aa6},\const{aSd}$ such that
\begin{subequations}
\begin{equation} \label{eq:m2pkwithGrad}
	\hat{c}\into (u+1)^{p+m_2+k-1} \leq \frac{p-1}{(p+m_1-1)^2}\into \abs*{\nabla (u+1)^{\frac{p+m_1-1}{2}}}^2+\const{a2} \quad
	\text{on }(0,\TM), \text{ provided \ref{M3}},
\end{equation}
\begin{equation} \label{eq:m3plwithGrad}
	\hat{c} \into (u+1)^{p+m_3+l-1} \leq \frac{p-1}{(p+m_1-1)^2} \into \abs*{\nabla (u+1)^{\frac{p+m_1-1}{2}}}^2+\const{a4} \quad
	\text{on }(0,\TM), \text{ provided \ref{M5}},
\end{equation}
\begin{equation} \label{eq:m3witheHpGN}
	  \hat{c} \into (u+1)^{p+m_3-1}\into (u+1)^l \leq \varepsilon_1\into (u+1)^{p+m_3+l-1}+ \frac{p-1}{(p+m_1-1)^2} \into \abs*{\nabla (u+1)^{\frac{p+m_1-1}{2}}}^2 + \const{a6} \quad\text{on }(0,\TM), 
	\end{equation}
	\begin{equation} \label{eq:m3pwwithHpGN}
	   \hat{c} \into (u+1)^{p+m_3-1}w \leq \varepsilon_2\into (u+1)^{p+m_3+l-1}+ \frac{p-1}{(p+m_1-1)^2}\into \abs*{\nabla (u+1)^{\frac{p+m_1-1}{2}}}^2 + \const{a4s} \quad\text{on }(0,\TM),
	\end{equation}
	\begin{equation} \label{eq:gradwithPhi}
		-\into  \abs*{\nabla (u+1)^{\frac{p+m_1-1}{2}}}^2 \leq \const{aa6}-\const{aSd}\tonda*{\into (u+1)^p}^{\frac{p+m_1-1}{p\theta_1}} \quad\text{on }(0,\TM).
	\end{equation}
\end{subequations}
\begin{proof} 
Let us show \eqref{eq:m2pkwithGrad}. Under assumption \ref{M3}, by
taking into account \eqref{eq:theta}, \eqref{eq:sigmatheta} and the boundedness of the mass (Lemma \ref{boundednessMass}), we can derive through the Gagliardo--Nirenberg and Young’s inequalities this bound
\begin{equation*}
\begin{split}
	\hat{c}\into (u+1)^{p+m_2+k-1}=&\hat{c}\norm*{(u+1)^\frac{p+m_1-1}{2}}_{L^{\frac{2(p+m_2+k-1)}{p+m_1-1}}(\Omega)}^\frac{2(p+m_2+k-1)}{p+m_1-1}\\
	\leq& \const{s2}\tonda*{\norm*{\nabla(u+1)^\frac{p+m_1-1}{2}}_{L^2(\Omega)}^{\theta}  \norm*{(u+1)^\frac{p+m_1-1}{2}}_{L^\frac{2}{p+m_1-1}(\Omega)}^{1-\theta}+\norm*{(u+1)^\frac{p+m_1-1}{2}}_{L^\frac{2}{p+m_1-1}(\Omega)}}^{\sigma}\\
	\leq& \const{s3}\tonda*{\into\abs*{\nabla (u+1)^\frac{p+m_1-1}{2}}^2}^\frac{\sigma\theta}{2}+\const{s4} \leq \frac{p-1}{(p+m_1-1)^2} \into\abs*{\nabla (u+1)^\frac{p+m_1-1}{2}}^2 + \const{a2}\quad \text{on }(0,\TM).
\end{split}
\end{equation*}
(We precise that we have made use of 
\begin{equation}\label{disab}
(A+B)^s\leq 2^s(A^s+B^s) \quad \text{for all } A,B,s>0,
\end{equation}
and we might employ this algebraic relation without mentioning it if not necessary.)  

Similarly to what we have done before, by supposing \ref{M5}, estimate \eqref{eq:m3plwithGrad} is obtained by applying in this case \eqref{eq:thetabar} and \eqref{eq:sigmathetabar}, so entailing 
\begin{equation*}
\begin{split}
	\hat{c}\into (u+1)^{p+m_3+l-1}=&\hat{c}\norm*{(u+1)^\frac{p+m_1-1}{2}}_{L^{\frac{2(p+m_3+l-1)}{p+m_1-1}}(\Omega)}^\frac{2(p+m_3+l-1)}{p+m_1-1} \\
	\leq& \const{S2}\tonda*{\norm*{\nabla(u+1)^\frac{p+m_1-1}{2}}_{L^2(\Omega)}^{\theta_1}  \norm*{(u+1)^\frac{p+m_1-1}{2}}_{L^\frac{2}{p+m_1-1}(\Omega)}^{1-\theta_1}+\norm*{(u+1)^\frac{p+m_1-1}{2}}_{L^\frac{2}{p+m_1-1}(\Omega)}}^{\sigma_1} \\
	\leq& \const{S3}\tonda*{\into\abs*{\nabla (u+1)^\frac{p+m_1-1}{2}}^2}^\frac{\sigma_1\theta_1}{2}+\const{S4} 
 \leq \frac{p-1}{(p+m_1-1)^2} \into\abs*{\nabla (u+1)^\frac{p+m_1-1}{2}}^2 + \const{a4}\quad \text{for all } t \in (0,\TM).
\end{split}
\end{equation*}
In order to prove \eqref{eq:m3witheHpGN} we distinguish the cases $l\leq1$ and $l>1$. If $0<l\leq1$ the boundedness of the mass, given by Lemma \ref{boundednessMass}, and Young's inequality yield
\begin{equation*}
\begin{split}
	\hat{c}\into (u+1)^{p+m_3-1}\into (u+1)^l &\leq \const{asd}\into (u+1)^{p+m_3-1} \leq \varepsilon_1\into (u+1)^{p+m_3+l-1}+\const{a6} \quad\text{on }(0,\TM).
\end{split}
\end{equation*}
On the other hand, if $l>1$ we first have to introduce $q>\max{\{l,m_3+l-1\}}$. Subsequently, by applying twice the H\"{o}lder inequality, we get by relying on Lemma \ref{lem:disab}  
\begin{equation}\label{DoppiaHolder}
\begin{split}
	\hat{c}\into (u+1)^{p+m_3-1}\into (u+1)^l &\leq \const{Gs3}\quadra*{\into(u+1)^{p+m_3+l-1} }^{\frac{p+m_3-1}{p+m_3+l-1}} \quadra*{\tonda*{\into (u+1)^q}^{\frac pq +1}}^{\frac{l}{p+q}}   \\
	& \leq \varepsilon_1 \into(u+1)^{p+m_3+l-1} + \varepsilon_1 \tonda*{\into (u+1)^q}^{\frac pq +1} +\const{as6} \quad\text{for all }t \in (0,\TM).
\end{split}
\end{equation}
Now we focus on the second integral of the right-hand side of \eqref{DoppiaHolder}. By exploiting \eqref{eq:thetatilde}, \eqref{eq:sigmathetatilde} and Lemma \ref{boundednessMass}, a combination of the Gagliardo–Nirenberg and Young’s inequalities gives 
\begin{equation*}
\begin{split}
	\varepsilon_1 \tonda*{\into (u+1)^q}^{\frac pq +1}=& \varepsilon_1\norm*{(u+1)^\frac{p+m_1-1}{2}}_{L^{\frac{2q}{p+m_1-1}}(\Omega)}^\frac{2(p+q)}{p+m_1-1}\\
	\leq& \const{ss2}\tonda*{\norm*{\nabla(u+1)^\frac{p+m_1-1}{2}}_{L^2(\Omega)}^{\theta_4}  \norm*{(u+1)^\frac{p+m_1-1}{2}}_{L^\frac{2}{p+m_1-1}(\Omega)}^{1-\theta_4}+\norm*{(u+1)^\frac{p+m_1-1}{2}}_{L^\frac{2}{p+m_1-1}(\Omega)}}^{\sigma_2}\\
	\leq& \const{ss3}\tonda*{\into\abs*{\nabla (u+1)^\frac{p+m_1-1}{2}}^2}^\frac{\sigma_2\theta_4}{2}+\const{ss4} \leq \frac{p-1}{(p+m_1-1)^2} \into\abs*{\nabla (u+1)^\frac{p+m_1-1}{2}}^2 + \const{a6}\quad \text{on }(0,\TM),
\end{split}
\end{equation*}
which plugged in the previous one concludes the proof.

As to \eqref{eq:m3pwwithHpGN}, again by means of the Young inequality and relation \eqref{ellipticreg}, recalling the definition of $g$ in \eqref{eq:disf} we can write   
\begin{equation}\label{EllRegSTim}
\begin{split}
	\hat{c}\into (u+1)^{p+m_3-1}w & \leq \rho \into(u+1)^{p+m_3+l-1} + \const{sgar} \into w^{\frac{p+m_3-1}{l}+1}\\ 
 &\leq \rho \into(u+1)^{p+m_3+l-1}+\rho \into g(u)^{\frac{p+m_3+l-1}{l}}+ C_\rho \tonda*{\into g(u)}^{\frac{p+m_3+l-1}{l}}\\
	&\leq \varepsilon_2\into(u+1)^{p+m_3+l-1}+\const{asf2}\tonda*{\into (u+1)^l}^{\frac{p+m_3+l-1}{l}} \quad \text{for all } t \in (0,\TM).
\end{split}
\end{equation}
For $l\leq 1$, due to Lemma \ref{boundednessMass} the last integral of the right-hand side of \eqref{EllRegSTim} is bounded, and we have 
\begin{equation*}
\hat{c}\into (u+1)^{p+m_3-1}w\leq \varepsilon_2\into(u+1)^{p+m_3+l-1} +\const{a4s} \quad \text{on }(0,\TM).
\end{equation*}
Instead, for $l>1$, by virtue of \eqref{eq:thetahat}, Lemma \ref{boundednessMass} and \eqref{eq:sigmathetahat}, the Gagliardo--Nirenberg and Young's inequalities imply 
\begin{equation*}
\begin{split}
\const{asf2} \tonda*{\into (u+1)^l}^{\frac{p+m_3+l-1}{l}}&= \const{asf2}\norm*{(u+1)^{\frac{p+m_1-1}{2}}}^{\sigma_1}_{L^{\frac{2l}{p+m_1-1}}(\Omega)}\\
&\leq \const{asf4}\tonda*{\norm*{\nabla(u+1)^{\frac{p+m_1-1}{2}}}^{\theta_2}_{L^{2}(\Omega)}\norm*{(u+1)^{\frac{p+m_1-1}{2}}}^{1-\theta_2}_{L^{\frac{2}{p+m_1-1}}(\Omega)}+\norm*{(u+1)^{\frac{p+m_1-1}{2}}}_{L^{\frac{2}{p+m_1-1}}(\Omega)}   }^{\sigma_1}\\
&\leq \const{asf5}\tonda*{\into\abs*{\nabla (u+1)^{\frac{p+m_1-1}{2}}}^2}^{\frac{\sigma_1\theta_2}{2}}+\const{asf6}\\
&\leq \frac{p-1}{(p+m_1-1)^2} \into \abs*{\nabla (u+1)^{\frac{p+m_1-1}{2}}}^2 +\const{a4s}  \quad \text{for all } t \in (0,\TM),
\end{split}
\end{equation*}
and we conclude by inserting this gained bound into \eqref{EllRegSTim}.

Finally, let us prove \eqref{eq:gradwithPhi}. By taking into account \eqref{eq:thetaunder} and Lemma \ref{boundednessMass}, a further application of the Gagliardo--Nirenberg inequality 
leads to
\begin{equation*}
\begin{split}
    \into (u+1)^p=&\norm*{(u+1)^\frac{p+m_1-1}{2}}_{L^{\frac{2p}{p+m_1-1}}(\Omega)}^{\frac{2p}{p+m_1-1}}\\ 
    \leq& \const{d_1}\tonda*{\norm*{\nabla(u+1)^\frac{p+m_1-1}{2}}_{L^2(\Omega)}^{\theta_3} \norm*{(u+1)^\frac{p+m_1-1}{2}}_{L^\frac{2}{p+m_1-1}(\Omega)}^{1-\theta_3}+\norm*{(u+1)^\frac{p+m_1-1}{2}}_{L^\frac{2}{p+m_1-1}(\Omega)}}^{\frac{2p}{p+m_1-1}}\\
	\leq& \const{d_3}\tonda*{\into\abs*{\nabla (u+1)^\frac{p+m_1-1}{2}}^2}^\frac{p\theta_3}{p+m_1-1}+\const{d_4}\leq \const{d5}\tonda*{\into\abs*{\nabla (u+1)^\frac{p+m_1-1}{2}}^2 +1}^\frac{p\theta_3}{p+m_1-1}  \quad \text{on }(0,\TM),
\end{split}
\end{equation*}
obtaining the claim after basic manipulations in the previous estimate. 
\end{proof}
\end{lemma}
\begin{lemma}\label{lem:validforboth}
For the value of $\Bar{p}>1$ found in Lemma \ref{lem:p}, let for $p>\Bar{p}$ define the functional $\varphi(t)$ by
\begin{equation*}
\varphi(t) := \frac1p \into (u+1)^p,
\end{equation*}
and for $j\in\{m_2,m_3\}$, let 
\begin{equation} \label{eq:funcor}
	F_j(u) := \int_0^u \hat{u} \tonda*{\hat{u}+1}^{p+j-3} d \hat{u}.
\end{equation} Then for all $t\in(0,\TM)$ and some $\const{aa3}$ and $\const{aa4}$
\begin{equation}\label{eq:validforboth}
	\varphi'(t)\leq -\frac{4(p-1)}{(p+m_1-1)^2}\into \abs*{\nabla (u+1)^{\frac{p+m_1-1}2}}^2 
	- (p-1)\chi \into F_{m_2}(u)  \Delta v
	+ (p-1)\xi \into F_{m_3}(u) \Delta w-\const{aa3}\into (u+1)^{p+r-1}+\const{aa4}.
\end{equation}
\begin{proof}
Let us differentiate the energy functional; we have
	\begin{align*}
	\varphi'(t) = &\into (u+1)^{p-1} u_t \notag \\
	           =& \into (u+1)^{p-1} \nabla \cdot \tonda*{(u+1)^{m_1-1}\nabla u}
	           -  \chi \into (u+1)^{p-1} \nabla \cdot \tonda*{u(u+1)^{m_2-1}\nabla v}  \notag \\
				& + \xi \into (u+1)^{p-1} \nabla \cdot \tonda*{u(u+1)^{m_3-1}\nabla w}
				+ \lambda \into (u+1)^{p-1} u - \mu \into (u+1)^{p-1} u^r \quad \text{on }(0,\TM). \notag \\
				\intertext{Due to the divergence theorem, the three first integrals on the right hand side make that the above identity becomes}
				\varphi'(t) =& -(p-1)\into (u+1)^{p+m_1-3} \abs*{\nabla u}^2  
				+ (p-1) \chi \into u(u+1)^{p+m_2-3} \nabla u \cdot \nabla v \notag \\
				&- (p-1) \xi \into u(u+1)^{p+m_3-3} \nabla u \cdot \nabla w
				+  \lambda \into (u+1)^{p-1} u - \mu \into (u+1)^{p-1} u^r  \quad \text{for all } t \in (0,\TM). \notag \\
				\intertext{By recalling \eqref{eq:funcor}, we rewrite the previous expression as}
				\varphi'(t)=&-\frac{4(p-1)}{(p+m_1-1)^2}\into \abs*{\nabla (u+1)^{\frac{p+m_1-1}2}}^2 
				+  (p-1) \chi \into \nabla F_{m_2}(u) \cdot \nabla v
				-  (p-1) \xi \into \nabla F_{m_3}(u) \cdot \nabla w \notag \\
				&+ \lambda \into (u+1)^{p-1} u - \mu \into (u+1)^{p-1} u^r \quad \text{on $(0,\TM)$},\notag \\
                 \intertext{and again the divergence theorem provides}
				\varphi'(t)=& -\frac{4(p-1)}{(p+m_1-1)^2}\into \abs*{\nabla (u+1)^{\frac{p+m_1-1}2}}^2 
				-(p-1)\chi \into F_{m_2}(u)  \Delta v
				+ (p-1)\xi \into F_{m_3}(u) \Delta w \\
				&+ \lambda \into (u+1)^{p-1} u - \mu \into (u+1)^{p-1} u^r, \quad \text{for all $t \in (0,\TM)$}.
\end{align*}
The inequality
\begin{equation*} 
	\frac{1}{2^r} \into (u+1)^{p-1+r} \leq \into (u+1)^{p-1} u^r + \into (u+1)^{p-1}\quad \text{on }(0,\TM),
\end{equation*}
justified by inequality \eqref{disab}, leads to the further bound
\begin{equation}\label{eq:disphibe}
\begin{split}
	\varphi'(t) \leq & -\frac{4(p-1)}{(p+m_1-1)^2}\into \abs*{\nabla (u+1)^{\frac{p+m_1-1}2}}^2 
	- (p-1)\chi \into F_{m_2}(u)  \Delta v
	+ (p-1)\xi \into F_{m_3}(u) \Delta w \\
	&+ \lambda \into (u+1)^p + \mu \into (u+1)^{p-1} - \frac{\mu}{2^r} \into (u+1)^{p-1+r} \quad \text{for all } t\in (0,\TM).
\end{split}     
\end{equation}
Finally, a double application of Young's inequality (recalling $r>1$) gives 
\[
\lambda \into (u+1)^p + \mu \into (u+1)^{p-1} \leq \frac{\mu}{2^{r+1}} \into (u+1)^{p-1+r} + \const{TTT} \quad \text{on } (0,\TM),
\]
which, in conjunction with inequality \eqref{eq:disphibe}, leads to \eqref{eq:validforboth}.
\end{proof}
\end{lemma}
From now on we will distinguish the analyses for the local problem and the nonlocal one.
\subsection{Study of the local problem \eqref{eq:localproblem}}
Let us separate the elliptic and parabolic cases, i.e., $\tau = 0$ and $\tau = 1$ respectively.

\subsubsection{The parabolic-elliptic case: \texorpdfstring{$\tau = 0$}{\textepsilon=0}}

\begin{lemma} \label{lem:localelimitatotau0}
For the value of $\Bar{p}>1$ found in Lemma \ref{lem:p}, let conditions \eqref{eq:disf}, \ref{M1}, \ref{M2} and \ref{M3} be satisfied. Then $u\in L^\infty((0,\TM),L^p(\Omega))$ for all $p>\Bar{p}$.
% \begin{equation*}
% 	\into u^p\leq C 
% \end{equation*}
\begin{proof}
The second and third equations in \eqref{eq:localproblem} allow the substitution of $\Delta v$ and $\Delta w$ in \eqref{eq:validforboth}, obtaining once nonpositive terms are dropped 
\begin{equation}\label{eq:logbound}
\begin{split}
\varphi'(t)\leq & -\frac{4(p-1)}{(p+m_1-1)^2}\into \abs*{\nabla (u+1)^{\frac{p+m_1-1}2}}^2 
	- (p-1)\chi \into F_{m_2}(u) \tonda*{\beta v - f(u)} 
	+ (p-1)\xi \into F_{m_3}(u) \tonda*{\delta w - g(u)} \\
    & \quad -\const{aa3}\into (u+1)^{p+r-1}+\const{aa4} \\
    \leq & -\frac{4(p-1)}{(p+m_1-1)^2} \into \abs*{\nabla (u+1)^{\frac{p+m_1-1}2}}^2 + (p-1)\chi \into F_{m_2}(u) f(u)
    +  \const{locrep1} \into F_{m_3}(u) w -  (p-1)\xi \into F_{m_3}(u) g(u)  \\
    & \quad - \const{aa3} \into (u+1)^{p+r-1} + \const{aa4} \quad \text{on } (0,\TM). 
\end{split}
\end{equation}
With some calculations, by using the definition in \eqref{eq:funcor}, we find that
			\begin{equation} \label{eq:CorFuncbounds}
				0 \leq \frac1{p+j-1} u^{p+j-1} \leq F_j(u) \leq \frac1{p+j-1} \tonda*{(u+1)^{p+j-1} - 1}.
			\end{equation}
			By considering \eqref{eq:disf}, \eqref{eq:logbound}, and \eqref{eq:CorFuncbounds}, we have
\begin{equation}
        \begin{split}
    			\varphi'(t)  \leq & -\frac{4(p-1)}{(p+m_1-1)^2}\into \abs*{\nabla (u+1)^{\frac{p+m_1-1}2}}^2 
    			+ \const{loccan} \into (u+1)^{p+m_2+k-1}+ \const{locrepn} \into (u+1)^{p+m_3-1} w \\
    			& - \const{locrepnn} \into u^{p+m_3 -1} (u+1)^l - \const{aa3} \into (u+1)^{p-1+r}+ \const{aa4}  \\		
    			\leq & -\frac{4(p-1)}{(p+m_1-1)^2}\into \abs*{\nabla (u+1)^{\frac{p+m_1-1}2}}^2 
    			+ \const{loccan} \into (u+1)^{p+m_2 + k-1} + \const{locrepn} \into (u+1)^{p+m_3-1} w  \\
    			& - \const{locrepnnn} \into (u+1)^{p+m_3 +l -1}  +  \const{locrepnn} \into (u+1)^l
    			- \const{aa3} \into (u+1)^{p-1+r}+ \const{aa4} \quad \text{for all } t \in (0,\TM), \label{eq:LocalbeforeGNYoungHoelder} 
    			\end{split}
\end{equation}
			where we have rearranged the term proportional to $-\into u^{p+m_3 -1} (u+1)^l$ through \eqref{disab}.
    Let us deal with the two terms $\into (u+1)^l$ and $\into (u+1)^{p+m_3-1} w$. For every $l>0$, we have thanks to the Young inequality that for all $\hat{c}>0$
   \begin{equation}\label{uL}
   \hat{c} \into (u+1)^l \leq \frac{\const{aa3}}{2} \into (u+1)^{p-1+r}+ \const{SSS} \quad \text{on } (0,\TM),    
   \end{equation}
   (for $l\leq 1$ the boundedness of the mass (see again Lemma \ref{boundednessMass}) would provide a sharper estimate), whereas through \eqref{eq:m3pwwithHpGN} applied with 
   $\varepsilon_2<\const{locrepnnn}$, 
   from \eqref{eq:LocalbeforeGNYoungHoelder} we derive for all $t \in (0,\TM)$
			\begin{equation}\label{LaStessaBis}
				\varphi'(t)  \leq -\const{S1} \into \abs*{\nabla (u+1)^{\frac{p+m_1-1}2}}^2 
				+\const{loccan} \into (u+1)^{p+m_2 + k-1}- \const{ab2} \into (u+1)^{p+m_3 +l -1} - \frac{\const{aa3}}{2} \into (u+1)^{p-1+r} + \const{ab4}.
			\end{equation}
			Now we use assumption \ref{M1} to write by Young's inequality
   \begin{equation*}
      \const{loccan} \into (u+1)^{p+m_2 + k-1} \leq \const{ab2} \into (u+1)^{p+m_3 +l -1} + \const{GGG} \quad \text{on } (0,\TM),     
   \end{equation*}
which, introduced into \eqref{LaStessaBis}, provides after neglecting the term from the logistic source 
 \begin{equation*}
				\varphi'(t)  \leq -\const{S1} \into \abs*{\nabla (u+1)^{\frac{p+m_1-1}2}}^2 
				+ \const{ab5} \quad \text{for all } t\in (0,\TM).
			\end{equation*}
			The same estimate, up to constants, is obtained with assumption \ref{M2} by controlling $\into (u+1)^{p+m_2 + k-1}$ with $\into (u+1)^{p-1+r}$, invoking the Young inequality, or relying on assumption \ref{M3} by exploiting
            \eqref{eq:m2pkwithGrad}.
			
			Finally, inequality \eqref{eq:gradwithPhi} yields the initial problem
			\begin{equation*}%\label{eq:sistemafinale}
				\begin{dcases}
					\varphi'(t) \leq \const{ac1} - \const{ac2} \varphi(t)^{\frac{p+m_1-1}{p\theta_1}} & \text{on }(0,\TM),\\
					\varphi(0)=\frac 1p \into (u_0+1)^p,
				\end{dcases}
			\end{equation*}
			and ODI comparisons arguments ensure $\varphi(t)\leq \max{\graffa*{\varphi(0),\tonda*{\frac{\const{ac1}}{\const{ac2}}}^{\frac{p\theta_1}{p+m_1-1}}}}$ for all $t\in(0,\TM)$.
		\end{proof}
	\end{lemma}

	\subsubsection{The parabolic-parabolic case: \texorpdfstring{$\tau = 1$}{\textepsilon=1}}
	
	\begin{lemma}\label{lem:localelimitatotau1}
For the value of $\Bar{p}>1$ found in Lemma \ref{lem:p}, let conditions \eqref{eq:disf}, \ref{M2}, \ref{M3}, \ref{M4} and \ref{M5} be valid. Then $u\in L^\infty((0,\TM),L^p(\Omega))$ for all $p>\Bar{p}$.
		% \begin{equation*}
		% 	\into u^p\leq C 
		% \end{equation*}
		\begin{proof}
  We start dealing with bounds for the terms $\displaystyle \into F_{m_2}(u) \Delta v$ and $\displaystyle \into F_{m_3}(u) \Delta w$. Relations in \eqref{eq:CorFuncbounds} and Young's inequality provide for 
  all $t\in (0,\TM)$
			\begin{equation} \label{eq:F2Lv}
				- (p-1)\chi \into F_{m_2}(u) \Delta v \leq (p-1)\chi \into F_{m_2}(u) \abs*{\Delta v} 
				\leq \const{locpar2} \into (u+1)^{p+m_2-1} \abs*{\Delta v} 
				\leq \into (u+1)^{p+m_2+k-1} + \const{d1} \into \abs*{\Delta v}^{\frac{p+m_2-1+k}{k}}
			\end{equation}
			and analogously
			\begin{equation} \label{eq:F3Lw}
				(p-1)\xi \into F_{m_3}(u) \Delta w \leq \into (u+1)^{p+m_3+l-1} + \const{d2}\into \abs*{\Delta w}^{\frac{p+m_3-1+l}{l}}.
			\end{equation}
			For $z=v$, $\eta=\beta$, $\psi=f(u)$ (and using \eqref{eq:disf}) and $q=\frac{p+m_2+k-1}{k}$, inequality \eqref{eq:tau1} allows us to write
			\begin{equation}\label{regparv}
				\int_0^t e^s \into \abs*{\Delta v}^{\frac{p+m_2+k-1}{k}}\,ds \leq C_P\left(1 + \int_0^t e^s \into(u+1)^{p+m_2+k-1}\,ds\right) \quad \text{on }(0,\TM),
			\end{equation}
			and, reasoning similarly for the third equation in model \eqref{eq:localproblem},
			\begin{equation}\label{regparw}
				\int_0^t e^s \into \abs*{\Delta w}^{\frac{p+m_3+l-1}{l}}\,ds \leq \tilde{C}_P\left(1 +\int_0^t e^s \into(u+1)^{p+m_3+l-1}\,ds\right) \quad \text{for all } t \in (0,\TM).
			\end{equation}
			By substituting \eqref{eq:F2Lv} and \eqref{eq:F3Lw} into \eqref{eq:validforboth}, we obtain
			\begin{equation}\label{dislemma1local}
				\begin{split}
					\varphi'(t) \leq & - \frac{4(p-1)}{(p+m_1-1)^2}\into\abs*{\nabla (u+1)^\frac{p+m_1-1}{2}}^2 + \into(u+1)^{p+m_2+k-1} + \const{d1}\into \abs{\Delta v}^\frac{p+m_2+k-1}{k} +\into(u+1)^{p+m_3+l-1} \\ 
					&+\const{d2}\into \abs{\Delta w}^\frac{p+m_3+l-1}{l}  - \const{aa3}\into (u+1)^{p+r-1}+\const{aa4} \quad \text{on }(0,\TM).
				\end{split}
			\end{equation}
   Now we add to both sides of \eqref{dislemma1local} the term $\varphi(t)$, and successively we multiply by $e^t$. Since $e^t \varphi'(t) + e^t \varphi(t) =
\frac{d}{dt}(e^t \varphi(t))$, an integration over $(0,t)$ provides for all $t \in (0,\TM)$
\begin{equation}\label{StimaPrec}
                    \begin{split}
					e^t \varphi(t) \leq& \varphi(0) +\int_0^t e^s \left(- \frac{4(p-1)}{(p+m_1-1)^2}\into\abs*{\nabla (u+1)^\frac{p+m_1-1}{2}}^2 + \into(u+1)^{p+m_2+k-1} + \const{d1}\into \abs{\Delta v}^\frac{p+m_2+k-1}{k}\right. \\ 
					& \left. +\into(u+1)^{p+m_3+l-1} +\const{d2}\into \abs{\Delta w}^\frac{p+m_3+l-1}{l} +\frac 1p \into (u+1)^p - \const{aa3}\into (u+1)^{p+r-1}+\const{aa4}\right) \, ds.  
				\end{split}
			\end{equation}
			Relations \eqref{regparv} and \eqref{regparw}, in conjunction with Young's inequality (recall $r>1$) entail on $(0,\TM)$ once inserted into \eqref{StimaPrec}
			\begin{equation}\label{IneqExp}
				\begin{split}
					\frac{e^t}{p}\into u^p \leq e^t \varphi(t) \leq& \varphi(0) +\int_0^t e^s \left( -\frac{4(p-1)}{(p+m_1-1)^2}\into\abs*{\nabla (u+1)^\frac{p+m_1-1}{2}}^2 + \const{aas}\into(u+1)^{p+m_2+k-1} \right. \\ 
					& \left. +\const{abs}\into(u+1)^{p+m_3+l-1} -\frac{\const{aa3}}{2}\into (u+1)^{p+r-1}+\const{aaS4} \right) \,ds.
				\end{split}
			\end{equation}
		At this stage, analogously to what we have done in the previous lemma, if restriction \ref{M2} is supported with one between \ref{M4} and \ref{M5}, Young's inequality and \eqref{eq:m3plwithGrad} make that bound \eqref{IneqExp} is turned into 
\begin{equation}\label{StimaFinale} 
				e^t \into u^p\leq \const{sd} + \const{a_6e}(e^t -1) \quad \textrm{for all }  t \in (0,\TM).    
    \end{equation}
The same procedure applies if we take into account \ref{M3}, implying \eqref{eq:m2pkwithGrad}, with either \ref{M4} (by means of Young's inequality) or \ref{M5}, ensuring \eqref{eq:m3plwithGrad}. 

  Henceforth, also in this case inequality \eqref{StimaFinale} is derived and as a consequence, in each of the above situations, we conclude that
  \begin{equation*}
   \into u^p\leq \const{sdf}\quad \textrm{on } (0,\TM).
  \end{equation*}
	\end{proof}
	\end{lemma}
 
\subsection{Study of the nonlocal problem \eqref{eq:nonlocalproblem}}

\begin{lemma}\label{lem:nonlocalelimitato}
Let the hypotheses of Lemma \ref{lem:localelimitatotau0} be valid. Then $u\in L^\infty((0,\TM),L^p(\Omega))$ for all $p>\Bar{p}$.

\begin{proof}
  Starting from relation \eqref{eq:validforboth}, the second and the third equations in model \eqref{eq:nonlocalproblem} lead to 
\begin{equation*} %\label{eq:nonLocalbeforeGNYoungHoelder}
   \begin{split}
       \varphi'(t) \leq &-\frac{4(p-1)}{(p+m_1-1)^2}\into \abs*{\nabla (u+1)^{\frac{p+m_1-1}2}}^2 
				+ (p-1)\chi \into F_{m_2}(u) f(u) 
				+ \frac{(p-1)\xi}{|\Omega|} \into F_{m_3}(u) \into g(u) \\
				& -   (p-1)\xi \into F_{m_3}(u) g(u) - \const{aa3} \into (u+1)^{p-1+r} + \const{aa4}\\ 
    \leq & -\frac{4(p-1)}{(p+m_1-1)^2}\into \abs*{\nabla (u+1)^{\frac{p+m_1-1}2}}^2 
				+ \const{nonloccan} \into (u+1)^{p+m_2+k -1}  
				+ \const{nonlocrepint1} \into (u+1)^{p+m_3 -1}\into (u+1)^l\\ 
    &- \const{Ale} \into u^{p+m_3 -1} (u+1)^l - \const{aa3} \into (u+1)^{p-1+r} + \const{aa4}
    \quad \text{for all } t \in (0,\TM),
\end{split}
\end{equation*}
where we have taken in mind the properties of $F_j$ and $f$ and $g$ fixed in \eqref{eq:disf} and \eqref{eq:CorFuncbounds}, and dropped nonpositive terms. In turn, thanks again to \eqref{disab}, we have
\begin{align*}
				\varphi'(t) &\leq -\frac{4(p-1)}{(p+m_1-1)^2}\into \abs*{\nabla (u+1)^{\frac{p+m_1-1}2}}^2 
				+ \const{nonloccan} \into (u+1)^{p+m_2 + k -1} 
				+ \const{nonlocrepint1} \into (u+1)^{p+m_3 -1} \into (u+1)^l \notag \\
				&\phantom{\leq } - \const{Raf} \into (u+1)^{p+m_3+l-1} 
+ \const{nonlocreptintex} \into (u+1)^l - \const{aa3} \into (u+1)^{p-1+r} + \const{aa4} \quad \text{on }(0,\TM).
\end{align*}
With the aid of estimate \eqref{uL}, by exploiting \eqref{eq:m3witheHpGN} with $\varepsilon_1< \const{Raf}$, we obtain on $(0,\TM)$
			\begin{equation*}
				\varphi'(t)\leq -\const{S3} \into \abs*{\nabla (u+1)^{\frac{p+m_1-1}2}}^2 
				+ \const{nonloccan} \into (u+1)^{p+m_2 + k -1} 
				- \const{asda1}\into (u+1)^{p+m_3+l-1} - \frac{\const{aa3}}{2} \into (u+1)^{p-1+r} + \const{aa41} . 
			\end{equation*}
   Since up to constants the previous estimate coincides with \eqref{LaStessaBis}, we can follow the same arguments of Lemma \ref{lem:localelimitatotau0} to conclude.
\end{proof}
\end{lemma} 
\section{Proof of Theorems \ref{theo:localtau0} and \ref{theo:localtau1}}
For $\tau=0$, we use respectively Lemma \ref{lem:localelimitatotau0} and Lemma \ref{lem:nonlocalelimitato}, supported by the extension criterion in Lemma \ref{lem:existence}; in this way Theorem \ref{theo:localtau0} is established.
Theorem \ref{theo:localtau1} is evidently obtained by invoking Lemma \ref{lem:localelimitatotau1}.

	\bigskip 
	
	\subsubsection*{Acknowledgements}
All the authors are partially supported by the research
project {\em Analysis of PDEs in connection with real
phenomena}, CUP F73C22001130007, funded by
\href{https://www.fondazionedisardegna.it/}{Fondazione di Sardegna},
annuity 2021.	
 AC and SF are members of the {\em Gruppo Nazionale per l’Analisi Matematica, la Probabilità e le loro Applicazioni} (GNAMPA) of the Istituto Nazionale di Alta Matematica (INdAM). 
 RDF is member of the {\em Gruppo Nazionale per il Calcolo Scientifico} (GNCS) of INdAM and acknowledges financial support by PNRR e.INS Ecosystem of Innovation for Next Generation Sardinia (CUP F53C22000430001, codice MUR ECS0000038).
%under the National Recovery and Resilience Plan (NRRP), Mission 4 Component 2 Investment 1.5 - Call for tender No.3277 published on December 30, 2021 by the Italian Ministry of University and Research (MUR) funded by the European Union – NextGenerationEU. Project Code ECS0000038 – Project Title eINS Ecosystem of Innovation for Next Generation Sardinia – CUP F53C22000430001- Grant Assignment Decree No. 1056 adopted on June 23, 2022 by the Italian Ministry of University and Research (MUR).
SF is also partially supported by: Research project of MIUR (Italian Ministry of Education, University and Research) Prin 2022 
{\em Nonlinear differential problems with applications to real phenomena} (Grant Number: 2022ZXZTN2).

%\bibliography{Bibliography}{}
%\bibliographystyle{abbrv}
%\bibliographystyle{alpha}

\end{document}